\documentclass[12pt]{amsart}
\usepackage{amsmath,amssymb,url}

\newtheorem{thm}{Theorem}[section]

\newtheorem{lemma}[thm]{Lemma}
\newtheorem{cor}[thm]{Corollary}

\theoremstyle{remark}
\newtheorem{rem}[thm]{Remark}

\newcommand{\E}{{\mathbb E}}
\newcommand{\R}{{\mathbb R}}
\newcommand{\C}{{\mathbb C}}

\newcommand{\modo}[1]{{\left|#1\right|}}
\newcommand{\normo}[1]{{\left\|#1\right\|}}
\newcommand{\smodo}[1]{{\mathopen|#1\mathclose|}}
\newcommand{\snormo}[1]{{\mathopen\|#1\mathclose\|}}

\DeclareMathOperator{\curl}{curl}
\DeclareMathOperator{\divergence}{div}
\newcommand{\BMO}{\mathrm{BMO}}

\begin{document}
\title[Regularity of Navier-Stokes]
{Conditions implying regularity of the three dimensional Navier-Stokes
equation}
\author{Stephen Montgomery-Smith}
\makeatletter
\address{Department of Mathematics,
University of Missouri,
Columbia, MO 65211}
\email{stephen@math.missouri.edu}
\urladdr{\url{http://www.math.missouri.edu/~stephen}}
\thanks{The author was
partially supported
by an NSF grant.}

\begin{abstract}
We obtain logarithmic improvements for 
conditions for regularity of the
Navier-Stokes equation, similar to those of Prodi-Serrin or
Beale-Kato-Majda.  
Some of the proofs make use
of a stochastic approach involving Feynman-Kac like inequalities.
As part of the our methods, we give a different approach to a priori estimates
of Foia\c s, Guillop\'e and Temam.
\end{abstract}

\maketitle

\section{Introduction}

The version of the three dimensional
Navier-Stokes equation we study is the
differential equation in
$u = u(t) = u(x,t)$, where $ t\ge 0$, and $x \in \R^3$:
$$ \frac{\partial u}{\partial t}
   = \Delta u - L\divergence(u \otimes u),
   \quad u(0) = u_0.$$
Here $L$ denotes the Leray projection.
We will not usually be working with classical solutions.  We define
$u(t)$, $0 \le t \le T$,
to be a solution of the Navier-Stokes equation if, whenever $u(t_0)$ is 
sufficiently regular for a mild solution
$$ u(t) = e^{(t-t_0)\Delta} u(t_0) - \int_{t_0}^t
   e^{(t-s)\Delta} L\divergence(u(s) \otimes u(s)) \, ds $$   
to exist for 
$t \in [t_0,t_0+\tau)$ for some $\tau>0$, 
then $u(t)$ is equal to that mild solution in $[t_0,t_0+\tau)$.

We also use other ways to describe the three dimensional
Navier-Stokes equation.  First, let us 
denote the vorticity by $w = w(t) = w(x,t) = \curl u$.
If $w$ is sufficiently smooth then
$$ \frac{\partial w}{\partial t}
   = \Delta w - u \cdot \nabla w + w \cdot \nabla u ,
   \quad w(0) = \curl u_0 .$$

Another description is given by the so called magnetization variable
\cite{chorin}, \cite{montgomery-smith-pokorny}.  Let $m = m(t) = m(x,t)$ 
be a vector field satisfying an equation
$$ \frac{\partial m}{\partial t}
   = \Delta m - u \cdot \nabla m - m \cdot (\nabla u)^T ,
   \quad m(0) = u_0 + \nabla q_0 $$
for some scalar field $q_0=q_0(x)$.  (Here the superscript $T$ denotes
the transpose.)
Then under sufficient smoothness
assumptions we have that $u$ is the Leray projection of $m$.

A famous open problem is to prove regularity of the Navier-Stokes
equation, that is, if the initial data $u_0$ is in $L_2$ and is
regular (which in this paper we define to mean that it is
in the Sobolev spaces $W^{n,q}$ for some $2\le q<\infty$ and all positive
integers $n$), then the solution $u(t)$ is regular for all
$t\ge0$.  Such regularity would also imply uniqueness of the solution
$u(t)$.  Currently only the existence of weak solutions is known.
Also, it is known that for each regular $u_0$ that there exists
$t_0>0$ such that $u(t)$ is regular for $0\le t \le t_0$.
We refer the reader to 
\cite{cannone}, \cite{constantin-foias}, \cite{doering-gibbons},
\cite{lemarie-rieusset}, \cite{temam}.

In studying this problem, various conditions that imply regularity
have been obtained.  For example, the 
Prodi-Serrin conditions (\cite{prodi}, \cite{serrin})
state that for some
$2 \le p < \infty$, $3<q\le\infty$ with
$\frac2p + \frac 3q \le 1$ that
$$ \int_0^T \snormo{u(t)}_q^p \, dt < \infty $$
for all $T>0$.
If $u$ is
a weak solution to the Navier-Stokes equation 
satisfying a Prodi-Serrin condition,
with regular initial
data $u_0$, then $u$ is regular (see \cite{sohr}).
(Recently Escauriaza, Seregin and Sver\'ak \cite{escauriaza et al}
showed that the condition when $q=3$ and $p=\infty$ is also sufficient.)
This is a long way from what is currently known for the so called Leray-Hopf 
weak
solutions:
$$ \int_0^T \snormo{u(t)}_q^p \, dt < \infty $$
for $\frac2p + \frac 3q \ge \frac32$, $2 \le q \le 6$.

Another condition is that of Beale, Kato and Majda
\cite{beale-kato-majda}.  They show that regularity
follows from the condition
$$ \int_0^T \snormo{w(t)}_\infty \, dt < \infty $$
for all $T > 0$.
(In fact they proved this for the Euler equation, but the
proof works also for the Navier-Stokes equation with only small modifications.)
This was strengthened by Kozono and Taniuchi 
\cite{kozono-taniuchi} to show that regularity follows
from the condition
$$ \int_0^T \snormo{\nabla u(t)}_\BMO \, dt \approx
   \int_0^T \snormo{w(t)}_\BMO \, dt < \infty $$
for all $T > 0$, where here $\BMO$ denotes the space of functions with bounded
mean oscillation.

The purpose of this paper is threefold.  First, we would like to provide some
logarithmic improvements to these conditions.  Secondly, we would like to
present a stochastic approach to the Navier-Stokes equation, obtaining
our conditions using Feynman-Kac like inequalities.  Thirdly, we would like
to present a different process for creating 
estimates of
Foia\c s, Guillop\'e and Temam.

To this end, the first result of this paper is the logarithmic improvement
to the Prodi-Serrin conditions.

\begin{thm} \label{main prodi-serrin}
Let 
$2 < p < \infty$, $3<q < \infty$ with
$\frac2p+\frac3q=1$.
If $u$ is a solution to the Navier-Stokes equation satisfying
$$ \int_0^T \frac{\snormo{u(t)}_q^p}{1+\log^+\snormo{u(t)}_q}
   \, dt < \infty $$
for some $T>0$, 
then $u(t)$ 
is regular for $0 < t\le T$.
\end{thm}

We first present a proof of this result (and indeed of
a slightly stronger result) that uses a standard approach.  
Then 
we present a stochastic approach to the Navier-Stokes equation.
This is a kind of Lagrangian coordinates approach to the
Navier-Stokes equation, but with a probabilistic twist in that we follow
the path of each particle with a stochastic perturbation.  A similar approach
was adopted by Busnello, Flandoli and Romito \cite{busnello et al}.

From this we obtain the following
Beale-Kato-Majda type condition.
For $ 1 \le q < \infty$, define the function on $[0,\infty)$
$$ \Phi_q(\lambda) = \left(\frac{e^\lambda-1}{e-1}\right)^q .$$
Define the
$\Phi_q$-Orlicz norm on
any space of measurable functions by the formula
$$ \snormo f_{\Phi_q} = 
   \inf\left\{\lambda>0:
   \int \Phi_q(\smodo{f(x)}/\lambda) \, dx \le 1 \right\} .$$
(Thus the triangle inequality is a consequence of the fact that
$\Phi_q$ is convex, see \cite{kras-rutickii}.)  

\begin{thm} \label{main beale-kato-majda}
Let 
$1<q < \infty$, $3<r<\infty$, and $T>0$.
Suppose that $u$ is a solution to the Navier-Stokes equation satisfying
\begin{enumerate}
\item for all $T_0\in(0,T)$
$$\displaystyle \int_{T_0}^T \snormo{\nabla u(t)}_{\Phi_q}
   \, dt < \infty ,$$ 
and
\item
either $q<3$, or $\snormo{u(t)}_r < \infty$ for almost every $t \in [0,T]$.
\end{enumerate}
Then $u(t)$ 
is regular for $0 < t\le T$.
\end{thm}

\noindent
Note that since $\snormo\cdot_{\Phi_{q_1}} \le c \snormo\cdot_{\Phi_{q_2}}$
for $q_1 > q_2$, we may assume without loss of generality that $q>3/2$.
Next, if $3/2<q<3$, since $\snormo\cdot_q \le (e-1) \snormo\cdot_{\Phi_q}$, by
the Sobolev inequality we see that the second hypothesis is automatically
satisfied with $r=3q/(3-q)$.  Also, this hypothesis is always
satisfied for Leray-Hopf weak solutions with $r=6$.

Next we demonstrate how to obtain Theorem~\ref{main prodi-serrin}
from Theorem~\ref{main beale-kato-majda} using the following
result.  If $u$ is a solution to the Navier-Stokes equation, we define the
sets
$$ A^{n,q}_{T_0,T_1}(\lambda) =
   \{ t\in[T_0,T_1] :
   \snormo{\nabla^n u(t)}_q \ge \lambda \} .$$

\begin{thm} \label{main foias-guillope-temam}
Given $3 < q_1 \le q_2 \le \infty$, and a non-negative integer $n$,
there exists
constants $c_1,c_2,c_3>0$ such that
if $u(t)$, $0 \le t \le T_2$ 
is a solution to the Navier-Stokes equation, and if $0\le T_1 \le T_2$,
then for all
$r \in (0,\sqrt{T_2-T_1})$ we have
$$ \smodo{A^{n,q_2}_{T_1+r^2,T_2}(c_1 r^{3/q_2-n-1})} 
   \le c_2 \smodo{A^{0,q_1}_{T_1,T_2}(c_3 r^{3/q_1-1})} .$$
\end{thm}

A similar result that one can obtain (but we do not prove here) is 
that for positive integers $n$ we have
$ \smodo{A^{n,2}_{T_1+r^2,T_2}(c_1 r^{1/2-n})} 
   \le c_2 \smodo{A^{1,2}_{T_1,T_2}(c_3 r^{-1/2})} $.

\begin{cor} 
\label{cor foias-guillope-temam}
Under the hypotheses of Theorem~\ref{main foias-guillope-temam}, there
exists a constant $c>0$ with the following properties.
If
$\Theta(\lambda)$ is a positive increasing function of $\lambda\ge 0$,
define 
$$ \kappa = 
   \int_0^\infty \min\{(c\lambda^{-2}-T_0)^+,T_1\} \, d\Theta(\lambda) .$$
Then
$$
   \int_{T_0}^{T_1} \Theta(\snormo{\nabla^n u(s)}_{q_2}^{1/(1+n-3/q_2)}) \, ds
   \le 
   c \kappa + c \int_0^{T_1} \Theta(c\snormo{u(s)}_{q_1}^{1/(1-3/q_1)}) \, ds .
$$
Similarly,
$$
   \int_{T_0}^{T_1} \Theta(\snormo{\nabla^n u(s)}_2^{1/(n-1/2)}) \, ds
   \le 
   c \kappa + c \int_0^{T_1} \Theta(c\snormo{\nabla u(s)}_2^2) \, ds .
$$
\end{cor}

Since the Leray-Hopf weak solution to the Navier-Stokes equation satisfies
$\int_0^T \snormo{\nabla u(t)}_2^2 \, dt < \infty$, one can quickly 
recover the results of 
Foia\c s, Guillop\'e and Temam
\cite{foias et al} that say that
$\int_0^T \snormo{\nabla^n u(t)}_2^{1/(n-1/2)} \, dt < \infty$.

\section{Theorem~\ref{main prodi-serrin}}
\label{simple}

The hypothesis of Theorem~\ref{main prodi-serrin} imply that, 
given $\epsilon \in (0,T)$, there exists
$T_0 \in (0,\epsilon)$ with $u(T_0) \in  L_q$.
Let $T^* > T_0$ be the first point of non-regularity for $u(t)$.
It is well known that 
in order to show that $T^* > T$,
it is sufficient to show an \emph{a priori} estimate, that is
$\sup_{T_0 \le t < \min\{T^*,T\}} \snormo{u(t)}_q < \infty$.
This is because it is then possible to extend the regularity beyond
$T^*$ if $T^* \le T$.
Without loss of generality, it is sufficient to consider the case
$T = T^*$ (so as to obtain a contradiction).

\begin{proof}[Proof of Theorem~\ref{main prodi-serrin}]
We allow all constants to implicitly
depend upon $p$ and $q$.
Let us define quantities
\begin{align*}
v &= u \smodo u^{q/2-1} ,\\
A &= \sum_{i,j=1}^3 \left(\smodo u^{q/2-1} 
     \frac{\partial u_i}{\partial x_j} \right)^2 ,\\
B &= \sum_{i,j=1}^3 \left(\smodo u^{q/2-3} u_i \sum_{k=1}^3 u_k 
     \frac{\partial u_k}{\partial x_j}\right)^2
\end{align*}
Note that
\begin{gather*}
   \smodo{\nabla v}^2 := \sum_{i,j=1}^3 \left(
   \frac{\partial v_i}{\partial x_j}\right)^2
   \approx A + B,\\
   \sum_{i,j=1}^3 \frac\partial{\partial x_j}\left( \smodo u^{q-2} u_i \right)
   \frac{\partial u_i}{\partial x_j}
   \approx A + B,\\
   \sum_{i,j=1}^3 \left(\frac\partial{\partial x_j} 
   \left(\smodo u^{q-2} u_i\right)\right)^2
   \le
   c \smodo u^{q-2} \smodo{\nabla v}^2 .
\end{gather*}
We start with the Navier-Stokes equation, take the inner product with 
$u \smodo{u}^{q-2}$, and integrate over $\R^3$ to obtain
$$
    \snormo u_q^{q-1} \frac{\partial}{\partial t} \snormo u_q
    =
    \int \smodo u^{q-2} u \cdot \Delta u \, dx
    -
    \int \smodo u^{q-2} u \cdot L \divergence (u \otimes u) \, dx .
$$
Integrating by parts, we see that
$$ \int \smodo u^{q-2} u \cdot \Delta u \, dx
   =
   - \int 
   \sum_{i,j=1}^3 \frac\partial{\partial x_j}\left( \smodo u^{q-2} u_i \right)
   \frac{\partial u_i}{\partial x_j}
   \, dx
   \approx
   - \snormo{\nabla v}_2^2 ,$$
and
\begin{align*} 
   \int \smodo u^{q-2} u \cdot L \divergence (u \otimes u) \, dx
   &=
   \int \sum_{i,j=1}^3
        \frac\partial{\partial x_j} \left(\smodo u^{q-2} u_i\right)
        [L (u_j u)]_i \, dx \\
   &\le
   c \snormo{\smodo{u}^{q/2-1}}_s \snormo{\nabla v}_2 \snormo{L(u \otimes u)}_r
\end{align*}
where $r = 1+q/2$ and $s = (2q+4)/(q-2)$.  Now the Leray projection is a bounded
operator on $L_r$, and hence 
$\snormo{L(u \otimes u)}_r \approx \snormo u_{2+q}^2$.  Also
$\snormo{\smodo{u}^{q/2-1}}_s \approx \snormo u_{2+q}^{q/2-1}$.  Hence
$$ \int \smodo u^{q-2} u \cdot L \divergence (u \otimes u) \, dx
   \le c \snormo u_{2+q}^{1+q/2} \snormo{\nabla v}_2 
   = c \snormo v_{2+4/q}^{1+2/q} \snormo{\nabla v}_2 .$$
From the Sobolev and interpolation inequalities
$$ \snormo v_{2+4/q} \le c \snormo{ \smodo{\nabla}^{3/(q+2)} v}_2
   \le c \snormo v_2^{(q-1)/(q+2)} \snormo{\nabla v}_2^{3/(q+2)} ,$$
and hence
$$ \int \smodo u^{q-2} u \cdot L \divergence (u \otimes u) \, dx
   \le c
   \snormo v_2^{1-1/q} \snormo{\nabla v}_2^{1+3/q} .$$
Now apply Young's inequality $ab \le ((q-3)a^{2q/(q-3)} + (q+3)b^{2q/(q+3)})/2q$
for $a,b \ge 0$, to obtain
$$
   \int \smodo u^{q-2} u \cdot L \divergence (u \otimes u) \, dx
   \le 
   c_1 \snormo{\nabla v}_2^2 + c_2 \snormo v_2^{2(q-1)/(q-3)} ,
$$
where $c_1$ may be made as small as required by making $c_2$ larger.  
Hence
$$ \snormo u_q^{q-1} \frac\partial{\partial t}\snormo u_q
   \le c \snormo v_2^{2(q-1)/(q-3)},$$
that is,
$$ \frac\partial{\partial t}\snormo u_q
   \le c \snormo u_q^{p+1} ,$$
and so
$$ \frac\partial{\partial t}
   \log(1+\log^+\snormo u_q)
   \le 
   \frac{c \snormo u_q^p}{1+\log^+\snormo u_q} .$$
Integrating, we see that for $T_0 \le t < T$
$$ \log(1+\log^+\snormo{u(t)}_q) \le
   \log(1+\log^+\snormo{u(T_0)}_q) + 
   c \int_{T_0}^T \frac{\snormo{u(s)}_q^p}{1+\log^+\snormo{u(s)}_q} \, ds ,$$
which provides a uniform bound for $\snormo{u(t)}_q$.
\end{proof}

\begin{rem}
Note that this proof can easily be adapted to show that a sufficient condition
for regularity is that
$$
   \int_0^T \frac{\snormo{u(s)}_q^p}{\Theta(\snormo{u(s)}_q)} 
   \, ds < \infty ,$$
where $\Theta$ is any increasing function for which 
$$ \int_1^\infty \frac1{x\Theta(x)} \, dx = \infty .$$
\end{rem}

\section{A Priori Estimates}

This section is devoted to the proof of Theorem~\ref{main foias-guillope-temam}
and Corollary~\ref{cor foias-guillope-temam}
The proof is very similar to the proof Scheffer's Theorem \cite{scheffer} 
that states that the Hausdorff dimension
of the set of $t$ for which the solution $u(t)$ is not regular is
$1/2$.
The main tool is
the following result is due to 
Gruji\'c and Kukavica \cite{grujic-kukavica} (see also 
\cite{lemarie-rieusset 2}).

\begin{thm} \label{space analytic}
There exist constants $a,c>0$ and a function
$T:(0,\infty) \to (0,\infty)$, with $T(\lambda) \to \infty$ as $\lambda\to 0$,
with the following properties.  If $u_0 \in L_q(\R^3)$, then there is
a solution $u(t)$ $(0 \le t \le T(\snormo{u_0}_q))$
to the Navier-Stokes equation, with $u(0) = u_0$, and
$u(x,t)$ is the restriction of an analytic function 
$u(x+iy,t) + iv(x+iy,t)$ in the region
$\{x+iy \in \C^3 : \smodo y \le a \sqrt t\}$, and
$\snormo{u(\cdot+iy,t) + i v(\cdot+iy,t)}_q \le c \snormo{u_0}_q$ for
$\smodo y \le a \sqrt t$.
\end{thm}

\begin{proof}[Proof of Theorem~\ref{main foias-guillope-temam}]
First let us show that
there exists
a constants $c_1,c_3,c_4>0$ such that 
if $u(t)$, $t_0-r^2 \le t \le t_0$ 
is a solution to the Navier-Stokes equation,
and 
$\smodo{A^{0,q_1}_{t_0-r^2,t_0}(c_3 r^{3/q_1-1})} < c_4 r^2$, then
$\snormo{\nabla^n u(t_0)}_{q_2} < c_1 r^{3/q_2-n-1}$.

To see this,
Let us first consider the case when $t_0 = 0$ and $r = 1$.
By hypothesis, we see that there exists
$t \in [-1,-1+c_4]$
with
$ \snormo{u(t)}_{q_1} < c_3$.  
By Theorem~\ref{space analytic} and the appropriate Cauchy integrals, 
if $c_4$
is small enough, then there exists a constant $c_7>0$ such that
$\snormo{\nabla^n u(0)}_{q_2} < c_1$.  

Now, by replacing 
$u(x,t)$
by
$r^{-1} u(r^{-1}x,r^{-2}(t-t_0))$, we can relax the restriction $r=1$ and
$t_0=0$, and we obtain the statement we asserted.

Next, given $\epsilon>0$,
it is trivial to find a finite collection $t_1,\dots,t_N$
in $A = A^{n,q_2}_{T_1+r^2,T_2}(c_1 r^{3/q_2-n-1})$
such that the sets $[t_n-r^2,t_n]$ are disjoint, but the sets
$[t_n-r^2-\epsilon,t_n+\epsilon]$ cover $A$.
By the above observation,
$\smodo{A^{0,q_1}_{t_0-r^2,t_0}(c_3 r^{3/q_1-1})} \ge c_4 r^2$.

Hence
$$
   \frac{r^2}{r^2+2\epsilon}\smodo{A}
   \le
   N r^2
   <
   c_4^{-1} 
   \sum_{n=1}^N 
   \smodo{A^{0,q_1}_{t_n-r^2,t_n}(c_3 r^{3/q_1-1})} 
   \le
   c_4^{-1} 
   \smodo{A^{0,q_1}_{T_1,T_2}(c_3 r^{3/q_1-1})} .
$$
Since $\epsilon$ is arbitrary, the result follows.
\end{proof}

\begin{proof}[Proof of Corollary~\ref{cor foias-guillope-temam}]
We only prove the first inequality.  By 
Theorem~\ref{main foias-guillope-temam}, there exist constants 
$c_1,c_2,c_3>0$ such
that
\begin{align*}
   \int_{T_0}^{T_1}&\Theta(\snormo{\nabla^n u(s)}_{q_2}^{1/(1+n-3/q_2)}) \, ds\\
   &=
   \int_0^\infty \smodo{\{s\in[T_0,T_1] : 
   \snormo{\nabla^n u(s)}_{q_2}^{1/(1+n-3/q_2)}>\lambda\}} \, d\Theta(\lambda)\\
   &\le
   c_1\kappa + 
   \int_0^\infty \smodo{\{s\in[c_2 \lambda^{-2},T_1] : 
   \snormo{\nabla^n u(s)}_{q_2}^{1/(1+n-3/q_2)}>\lambda\}} \, d\Theta(\lambda)\\
   &\le
   c_1 \kappa +
   c_1 \int_0^\infty \smodo{\{s\in[0,T_1] : 
   \snormo{u(s)}_{q_1}^{1/(1-3/q_1)} > c_3 \lambda\}} \, d\Theta(\lambda)\\
   &= 
   c_1 \kappa + c_1 \int_0^{T_1} 
   \Theta(c_3^{-1}\snormo{u(s)}_{q_1}^{1/(1-3/q_1)}) \, ds .
\end{align*}
\end{proof}

\section{A Stochastic Description}

Let us give a little motivation.  Suppose that 
we defined $\varphi_{t_0,t_1}(x)$
to be $X(t_0)$, where $X$ satisfies
the equation
$$ dX(t) = u(X(t),t) \, dt,
   \qquad
   X(t_1) = x ,$$
then $\varphi_{t_0,t_1}$ would be the ``back to coordinates map'' that
takes a point at $t=t_1$ to where it was carried from by the flow of the
fluid at time $t=t_0$.  For the Euler equation, this provides a very effective
way to describe the solution, for example, the
equation for vorticity can be rewritten in a Lagrangian form:
$$ w(x,t) = w(\varphi_{0,t}(x),0)
   + \int_0^t w(\varphi_{s,t}(x),s) \cdot \nabla u(\varphi_{s,t}(x),s) \, ds .$$
Similarly, for the magnetization variable we have
$$ m(x,t) = m(\varphi_{0,t}(x),0)
   - \int_0^t m(\varphi_{s,t}(x),s) \cdot 
   (\nabla u(\varphi_{s,t}(x),s))^T \, ds .$$
For the Navier-Stokes equation this formula is not true, and the Laplacian term
can make things complicated.  One approach to dealing with this is described
in the paper by Constantin \cite{constantin}.  
However, we take a different approach using Brownian motion, using a
kind of
``randomly perturbed back to coordinates map.''
Such a method was already discussed in the paper
\cite{montgomery-smith-pokorny}, here we make the discussion more rigorous.
The author recently found out that a similar approach was followed by
Busnello, Flandoli and Romito in \cite{busnello et al}.

The hypothesis of Theorem~\ref{main beale-kato-majda} imply that, 
given $\epsilon \in (0,T)$, there exists
$t' \in (0,\epsilon)$ with $u(t') \in  L_r$.
Then by known results (for example Theorem~\ref{space analytic}), 
it follows
that there exists $0 < T_0<\epsilon$ such that 
$u(T_0) \in W^{n,r'}$ for all $r' \in [r,\infty]$ and positive integers $n$.
Furthermore, arguing as in Section~\ref{simple}, we only need to prove 
$\sup_{T_0 \le t < \min\{T^*,T\}} \snormo{u(t)}_r < \infty$
under the \emph{a priori} assumption that
the solution is regular for $t \in [T_0,T]$.

If $f\colon\R^3 \to \R$ is regular, and
$T_0 \le t_0 \le t_1 < T$,
define $A_{t_0,t_1} f(x) = \alpha(x,t_1)$, where 
$\alpha$ satisfies the transport equation
$$ \frac{\partial\alpha}{\partial t} = \Delta \alpha - u\cdot\nabla \alpha,
   \qquad
   \alpha(x,t_0) = f(x) .$$
Since $\divergence(u) = 0$, an easy integration by parts argument shows
that
$$ \frac\partial{\partial t} \int \alpha(x,t) \, dx = 0 ,$$
and hence if $f$ is also in $L_1$, then
$$ \int A_{t_0,t_1} f(x) \, dx = \int f(x) \, dx .$$
Since stochastic
differential equations traditionally move forwards in time, it will be 
convenient to consider a time reversed equation.
Let $b(t)$ be three dimensional Brownian motion.
For $T_0 \le t_0 \le t_1 < T_1$, define the random function
$\varphi_{t_0,t_1}\colon\R^3\to\R^3$ by
$\varphi_{t_0,t_1}(x) = X(-t_0)$, where $X$ satisfies the 
stochastic differential equation:
$$ dX(t) = -u(X(t),t) \, dt + \sqrt2 \, db(t),
   \qquad
   X(-t_1) = x .$$
It follows by the Ito Calculus \cite{karatzas-shreve} that
if $T_0 \le t_0 \le t_1 < T$, then
$$ A_{t_0,t_1} f(x) = \E f(\varphi_{t_0,t_1}(x)) .$$
(Here as in the rest of the paper, $\E$ denotes expected value.)
Note that if $f$ is also in $L_1$, then
$$ \int \E f(\varphi_{t_0,t_1}(x)) \, dx = \int f(x) \, dx .$$
Applying the usual dominated and monotone convergence theorems, it
quickly follows that the last equality is also true if $f$ is any
function in $L_1$, or if $f$
is any positive function.

Now let us develop the equations for the magnetization variable.  (The same
approach will also work for the vorticity.)
If we set $m(T_0) = u(T_0)$, then
we note that $m$ is the unique solution to the integral
equation
$$ m(t) = A_{T_0,t} u(T_0) - 
   \int_{T_0}^t A_{s,t} (m(s) \cdot (\nabla u(s))^T) \, ds 
   \quad (T_0 \le t < T).$$
Uniqueness follows quickly by the usual fixed point argument
over short intervals, 
remembering that $u(t)$ is regular for $T_0 \le t < T$.

Consider also the random quantity
$\tilde m = \tilde m(x,t)$ as the solution to the integral equation
for $T_0 \le t < T$
$$ \tilde m(x,t) = u(\varphi_{T_0,t}(x),T_0) -
   \int_{T_0}^t \tilde 
   m(\varphi_{s,t}(x),s) \cdot (\nabla u(\varphi_{s,t}(x),s))^T \, ds .$$
Again, 
it is very easy to show that a solution exists by using a fixed point
argument over short time intervals.
It is seen that $\E\tilde m$ satisfies the same equation as $m$, and
hence $\E \tilde m = m$.

Next, 
$\varphi_{t_0,t_1}(\varphi_{t_1,t_2}(x)) = \varphi_{t_0,t_2}(x)$,
since both are $Y(t_0)$ where $Y(t)$ is the solution to the integral equation
$$ Y(t) = \varphi_{t_1,t_2}(x) + \int_{t_1}^{t} u(Y(s),s)\, ds
       + \sqrt 2(b_{-t}-b_{-t_1}) .$$
Hence
$$ \tilde m(\varphi_{s_1,t}(x),s_1) - \tilde m(\varphi_{s_2,t}(x),s_2)
   =
   \int_{s_1}^{s_2} \tilde 
   m(\varphi_{s,t}(x),s) \cdot (\nabla u(\varphi_{s,t}(x),s))^T \, ds .$$
Thus,
by Gronwall's inequality, if $T_0 \le t < T$
$$ \smodo{\tilde m(x,t)}
   \le
   \exp\left(\int_{T_0}^t \smodo{\nabla u(\varphi_{s,t}(x),s)} \, ds\right)
   \smodo{u(\varphi_{T_0,t}(x),T_0)} .$$
(This is essentially the Feynman-Kac formula.)
The goal, then, is to find uniform estimates on the quantity
$$ \exp\left(\int_{T_0}^t \smodo{\nabla u(\varphi_{s,t}(x),s)} \, ds\right) .$$
This we proceed to do in the next section.

\section{Theorem~\ref{main beale-kato-majda}}

Let us fix $q$ and $r$ satisfying the hypothesis of 
Theorem~\ref{main beale-kato-majda}, 
and allow all constants to implicitly depend upon $q$ and $r$.
We retain the notation from the previous section, in particular the 
definitions of $T_0$, $T^*$ and $T$.

\begin{proof}[Proof of Theorem~\ref{main beale-kato-majda}]
Since $\snormo{u(t)}_r < \infty$ for almost every $t \in [0,T]$, by 
Theorem~\ref{main foias-guillope-temam},
we see that $\snormo{\nabla u(t)}_\infty < \infty$ 
for almost every $t \in [0,T]$.
Hence, there exists $\lambda>T_0^{-1}$ such
that
$$ \int_{B}
   \snormo{\nabla u(t)}_{\Phi_q} \, dt \le \frac1q ,$$
where 
$B = \{t \in [T_0,T] \colon \snormo{\nabla u(t)}_\infty \ge c_2 \lambda\}$.
Thus for $T_0 \le t < T$, we have that $\smodo{\tilde m(x,t)}$ is bounded
by
$$ e^{c_2 \lambda (t-T_0)}
   \exp\left(\int_{B\cap[T_0,t]} 
   \smodo{\nabla u(\varphi_{s,t}(x),s)} \, ds\right)
   \smodo{u(\varphi_{T_0,t}(x),T_0)} .$$
Hence by Jensen's and
H\"older's inequalities, $\snormo{m(t)}_r^r \le 
\int \E\smodo{\tilde m(t)}^r \, dx \le 
e^{c_2 q \lambda (t-T_0)}(N_r^r + N_{rq'}^r \tilde N^r)$, where
$q' = q/(q-1)$, 
$$
   N_s
   =
   \left(\int \E \smodo{u(\varphi_{T_0,t}(x),T_0)}^s\,dx\right)^{1/s} 
   = \snormo{u(T_0)}_s, $$
and
$$
   \tilde N^q
   =
   \int \E \left(
   \exp\left(q\int_{B \cap [T_0,t]}
   \modo{\nabla u(\varphi_{s,t}(x),s)} \, ds \right) - 1 \right)^q
   \, dx .
$$
Since the Orlicz norm satisfies the triangle inequality, we have
$$ \normo{\int_{B\cap[T_0,t]} \smodo{\nabla u(\varphi_{s,t}(\cdot),s)} \, ds}
    _{\Phi_q} \le \frac1q ,$$
that is, $\tilde N \le e-1 $.
Since $a^r + b^r \le (a+b)^r$ for $a,b \ge 0$,
we conclude that 
$$ \snormo{m(t)}_r \le \snormo{u(T_0)}_r + 
   (e-1) e^{c_2 \lambda (t-T_0)} \snormo{u(T_0)}_{rq'} .$$
As the Leray projection is a bounded operator on $L_r$ for $1<r<\infty$,
it follows that $\snormo{u(t)}_r$ is also uniformly bounded,
and the result follows.
\end{proof}

A second proof of Theorem~\ref{main prodi-serrin} now follows from
this next result.

\begin{lemma}  \label{upper bound for L_Phi}
There is a constant $c>0$ such that if $f$ is a measurable
function, then
$$
  \snormo f_{\Phi_q} \le c
  \left(\snormo f_q + 
  \frac{\snormo f_\infty}{1+\Phi_q^{-1}(({\snormo f_\infty}/{\snormo f_q})^q)}
  \right) .$$
\end{lemma}

\begin{proof}
Let us assume that $\snormo f_\infty = 1$, 
and set $a = \snormo f_q$, $b = \Phi_q^{-1}(a^{-q})$ and
$n = a+1/(1+b)$.  
Let $f^*:[0,\infty]\to[0,\infty]$ be the non-increasing rearrangement
of $\modo{f}$, that is,
$$ f^*(t) =  
   \sup \{ \lambda>0 : \smodo{\{ x:\smodo{f(x)}>\lambda\}} > t \} ,$$
so $\int F(\modo{f(x)}) \, dx = \int_0^\infty F(f^*(t)) \, dt$ for any
Borel measurable function $F$.  Notice that $f^*(t) \le \min\{1,a t^{-1/q}\}$.

Let us first
consider the case $a \le 1$, so that $b \ge 1$, $2n \ge 1/b$, and
$n \ge a$.
Then
$$
   \int \Phi_q(\modo{f(x)}/2n) \, dx
   \le
   \int_0^\infty \Phi_q(f^*(t)/2n) \, dt .$$
We split this integral up into three pieces.  First,
$$ 
   \int_0^{a^q} \Phi_q(f^*(t)/2n) \, dt 
   \le
   \int_0^{a^q} \Phi_q(b) \, dt 
   = 1.
$$
Next, since $(\Phi_q(\lambda))^{1/2q}$ is convex for $\lambda \ge 1$,
\begin{align*}
   \int_{a^q}^{a^q b^q} \Phi_q(f^*(t)/2n) \, dt 
   &\le
   \int_{a^q}^{a^q b^q} \Phi_q(abt^{-1/q}) \, dt \\
   &\le 
   \int_{a^q}^{a^q b^q} \frac{a^{2q}\Phi_q(b)}{t^2} \, dt  \\
   &\le 1.
\end{align*}
Next, for $t \ge a^q b^q$, $f^*(t) \le 1/b \le 2n$, and
$\Phi_q(\lambda) \le \lambda^q$ for $0 \le \lambda \le 1$, so
$$
   \int_{a^q b^q}^\infty \Phi_q(f^*(t)/2n) \, dt 
   \le
   \int_{a^q b^q}^\infty (f^*(t)/2n)^q \, dt 
   \le 1.
$$
Since $\Phi_q(\lambda/3) \le \Phi_q(\lambda)/3$ for
$\lambda\ge0$,
$$
   \int \Phi_q(\modo{f(x)}/6n) \, dx
   \le 1 ,
$$
that is, $\snormo f_{\Phi_q} \le 6 n$.

The case $a \ge 1$ (so $b \le 1$ and $2n \ge 1+2a$) is simpler, as
it is easy to estimate
$$
   \int_0^\infty \Phi_q(f^*(t)/2n) \, dt 
   \le
   \int_0^1 \Phi_q(1) \, dt 
   +
   \int_1^\infty (f^*(t)/2n)^q \, dt 
   \le 2.
$$

\end{proof}

\begin{proof}[Second proof of Theorem~\ref{main prodi-serrin}]
Applying Corollary~\ref{cor foias-guillope-temam} using the function
$$ \Theta(\lambda) = \frac{\lambda^2}{1+\log^+\lambda} ,$$
we obtain for all $T_0 \in (0,T)$
$$ \int_{T_0}^T \frac{\snormo{\nabla u(s)}_\infty}
   {1+\log^+\snormo{\nabla u(s)}_\infty} \, ds < \infty $$
and
$$ \int_{T_0}^T \frac{\snormo{\nabla u(s)}_q^{2q/(2q-3)}}
   {1+\log^+\snormo{\nabla u(s)}_q} \, ds < \infty .$$
Hence if $1<\alpha<2q/(2q-3)$ we have that
$$ \int_{T_0}^T \snormo{\nabla u(s)}_q^\alpha \, ds < \infty .$$
Next, considering the cases 
$\snormo{f}_\infty > \snormo{f}_q^\alpha$
and
$\snormo{f}_\infty \le \snormo{f}_q^\alpha$,
we see that
$$ \frac{\snormo f_\infty}{1+\Phi_q^{-1}(({\snormo f_\infty}/{\snormo f_q})^q)}
   \le c \left( 
   \snormo{f}_q^\alpha +
   \frac{\snormo{f}_\infty}
   {1+\log^+\snormo{f}_\infty} \right) .$$
Applying Lemma~\ref{upper bound for L_Phi}, we see
that the hypothesis of Theorem~\ref{main prodi-serrin} implies the
hypotheses of Theorem~\ref{main beale-kato-majda} with $q=r$.
\end{proof}

\section*{Acknowledgments}

The author wishes to extend his sincere gratitude to Michael Taksar for
help with understanding stochastic processes, and also to Pierre-Gilles
Lemari\'e-Rieusset for very helpful email discussions.

\end{document}